\documentclass{article}
\def\citep#1{\cite{#1}}

\def\be{\begin{equation}}
\def\ee{\end{equation}}

\begin{document}

\title{Bio-Inspired Computation: Success and Challenges of IJBIC}

\author{Xin-She Yang \\
School of Science and Technology,\\ Middlesex University, London NW4 4BT, United Kingdom.
\and 
Zhihua Cui \\
Complex System and Computational Intelligent Laboratory, \\ Taiyuan University of
Science and Technology, Taiyuan, P.R.China}

\date{}
\maketitle

\begin{abstract}
It is now five years since the launch of the International Journal of Bio-Inspired Computation (IJBIC).
At the same time, significant new progress has been made in the area of bio-inspired computation.
This review paper summarizes the success and achievements of IJBIC in the past five years,
and also highlights the challenges and key issues for further research.

\end{abstract}

{\bf Citation detail:}
X. S. Yang, and Z. H. Cui, Bio-inspired Computation: Success and Challenges of IJBIC,
{\it Int. J. Bio-Inspired Computation}, Vol.~6, No.~1, pp.1-6 (2014).

%% Begin of Main Text %%

\section{Introduction}

Bio-inspired computation, especially those based on swarm intelligence, has been in rapid
developments in the last two decades \citep{Ken,Yang,Cui,Cui2010,Yangbook2013,Gandomi2013CPSO}.
The literature has been expanding significantly
with diverse applications in almost all area of science, engineering and industries \citep{YangKoziel,Yangbook2013}.

The launch of the IJBIC in 2009 was a significant step
to provide a timely platform for researchers to
exchange ideas in bio-inspired computation. Initially, it aimed to publish quarterly, 4 issues per year.
One year later in 2010, the journal had to expand to 6 issues a year, due to the overwhelming
number of submissions and interests generated in the research communities. Since then, it becomes a
steady inflow of high quality submissions with bimonthly publications. In 2012, IJBIC was included
by Thomson Reuters in their Science Citation Index (SCI) Expanded, and the first impact factor 1.35 was announced earlier in 2013.
In the short 5 years, IJBIC has achieved its goals successfully.

Many other journals such as Soft Computing, Swarm Intelligence,
IEEE Transaction on Evolutionary Computation, Neural Computing and Applications, and Heuristics
all can include bio-inspired computation as a topic; however,
IJBIC was the only journal that specifically focused on bio-inspired computation.
It is no surprise that IJBIC has become successfully in such a short period.

Bio-inspired computation (BIC) belongs to the more general nature-inspired computation \citep{Yang},
and BIC includes swarm intelligence (SI) as its subset. The algorithms that are based on swarm intelligence
are among the most popular algorithms for optimisation, computational intelligence and data ming \citep{Yangbook2013}.
These SI-based algorithms include ant colony optimisation, bee algorithms such as artificial bee colony,
bat algorithm, cuckoo search, firefly algorithm, particle swarm optimization and others. In fact,
significant progress has been made in the last twenties years, and the publications in IJBIC just
provides a good snapshot of the latest developments.

The aim of this paper is two folds: to summarize the successful developments
in IJBIC and to highlight the trends and key issues concerning bio-inspired computation.
Therefore, this paper is organized as follows. Section 2 summarizes the latest developments
concerning bio-inspired computation, while Section 3 discusses the current trends.
Section 4 highlights the challenges and key issues in bio-inspired computation.
Finally, Section 5 concludes briefly with some topics for further research.

\section{Recent Developments in Bio-Inspired Computation}

The inaugural issue kicks off with the first paper by \cite{GrossDorigo}
about swarms of robots for group transport, followed by studies on particle swarm
and other methods.  Then, \cite{Deb2009} discussed the alternatives
to machine vision. Since then, the scopes and aims of IJBIC has expanded significantly,
and the studies published in the journal can be summarized into two main categories:
algorithm developments and applications.

\subsection{Algorithm Developments}

In terms of algorithm developments, there are
three strands: the developments of new algorithms,
improvements of existing algorithms, and the analysis of algorithms. Due to the complexity
of the problems in bio-inspired computation, it is not practical to systematically
analyze all relevant research in detail. In stead, we will briefly touch
relevant studies and highlight the most relevant progress in all three areas.
In addition, it is not possible to include all the papers published in the last five years
in IJBIC, and thus we only sample a fraction (e.g., about a quarter or more research papers)
so as to provide a timely snapshot of the diversity and depth of the current research and developments.

Among the development of new algorithms published in IJBIC,
\cite{Yang2010FA} introduced the new firefly algorithm and some stochastic functions,
which is the most cited paper of IJBIC with a total number of 151 citations in last
3 years since its publication in 2010. Later \cite{Yang2011BA} extended the bat algorithm for single objective optimisation to solve
multiobjective optimisation problems. In fact,
the firefly algorithm represents a successful
algorithm with diverse applications \citep{Fister,Yang2013FA,Fister2013b,SrivSoft,Gandomi2013FA}.

In addition, \cite{Adamatzky} introduced a new idea using slime mould to
compute planar shapes, which mimics important features of shortest paths in
route systems such as highways.
Furthermore, \cite{Shah} presented an intelligent water drops (IWD) algorithm,
though strictly speaking, IDW is not a bio-inspired algorithm, however, it does has some
characteristics of swarm intelligence.

On the other hand, \cite{Ray2011} used DNA computing to carry similarity-based fuzzy reasoning, while \cite{Saha} introduced the bacteria foraging algorithm to solve optimal FIR filter design problems. In addition, \cite{Cui2013} proposed the artificial plant optimisation algorithm with detailed case studies, especially
\cite{Cai2012} carried light responsive curve selection for photosynthesis operator
of artificial plant optimisation algorithm.

In addition to algorithm developments, more research activities have focused on the
improvements of existing algorithms. For example, \cite{Neri} introduced a variant of differential evolution for optimisation
in noisy environments, while \cite{Pandi} introduced the modified harmony search
for solving dynamic economic load dispatch problems. At the same time,
\cite{Xie2010} provided a brief survey
of artificial physics optimisation, while \cite{Roy2010} used differential harmony search algorithm to design a fractional-order PID controller.

Furthermore, \cite{Pal2011} used a modified invasive weed optimisation algorithm to carry out linear antenna array synthesis, while \cite{Layeb} improved the cuckoo search algorithm by introducing a novel quantum inspired cuckoo search
to solve knapsack problems. On the other hand, \cite{YangES} combined the eagle strategy with differential evolution to produce a new hybrid, and
\cite{Jamil} extended cuckoo search to study multimodal function optimisation.

Theoretical analysis in contrast is relatively weak in IJBIC, which is also true for
most other journals in the whole area. However, some papers can be classified
or partly classified as theoretical analyses because they intended to address the theoretical
and mathematical aspects of recent algorithms. For example,  \cite{Yang2011rev} provide an analysis about random walks and its role in nature-inspired metaheuristic
algorithms, while \cite{Xiao} discussed relationships between swarm intelligence and
artificial immune system. In addition, \cite{Green} discussed the elements of a network theory of adaptive complex systems, while \cite{Komatsu} investigated the dynamic diffusion in evolutionary optimised networks.

The lack of theoretical studies suggests that it is highly needed to address many problems
in the near future, and we hope that this paper can inspire more research in this area.

\subsection{Applications}

Applications are very diverse and form by far the vast majority of the
research papers published in IJBIC. In fact, more than two thirds of all the studies
are about applications of bio-inspired algorithms. Such diversity is one of key reasons that  drive the success
of the journal and may inspire more applications. The expanding literature makes
it unrealistic to analyze or even list all the papers in applications published in
the last five years. Therefore, we only select a fraction of the papers and try to
provide a representative subset of the current activities in bio-inspired computation.

For example,
\cite{Arumugam} carried out the optimal control of the steel annealing processes
by PSO, while \cite{Sivan} presented a study
in dynamic task scheduling with load balancing using
parallel orthogonal particle swarm optimisation.
In addition, \cite{Khabb} used the imperialist competitive algorithm
for minimising bit error rates, and \cite{Singh} studied the damping of low frequency oscillations in power system network using swarm intelligence tuned a fuzzy controller,
while \cite{ZhangCai} used a hybrid PSO to study  economic dispatch problems.

Furthermore,
\cite{Yampol} applied bio-inspired algorithm to the problem of integer factorisation, while \cite{Sheta} designed analogue filters using differential evolution.
Then, \cite{Laalaoui} used an ACO approach with learning to solve preemptive scheduling tasks.

On the other hand, \cite{Bonyadi} used a dynamic max-min ant system to solve the travelling salesman problem, and \cite{Bessedik} investigated graph colouring using bees-based algorithm.

A few special issues have been edited to address a subset of active and special research topics. For example, \cite{Tanimoto} edited a special issue on evolving world: games, complex networks, and agent simulations. In addition, \cite{Liu2012} edited a special issue on computational intelligence and its applications in engineering, and \cite{Natar} carried out a comparative study of cuckoo search
and bat algorithm for bloom filter optimisation in spam filtering. In terms of the summaries of more recent developments, \cite{Yang2012SI} edited a special issue on metaheuristics and swarm intelligence in engineering and
industry, where \cite{Marichel} presented an improved hybrid cuckoo search for solving
permutation flow shop scheduling problems.

The diversity of the applications can be seen from a wide range of topics concerning applications. \cite{Handa} used neuroevolution with manifold learning for playing Mario, while \cite{Srivastava} solved test sequence optimisation using cuckoo search, and later \cite{Srivastava2012b} used an approach based on cuckoo search to estimate software test effort. Furthermore, \cite{Sing2013} used artificial bee colony algorithm to optimise coordination
of electro-mechanical-based overcurrent relays. In the area of image processing, \cite{Dey} used cuckoo search to optimise the scaling factors in electrocardiogram signal watermarking.

In the area of web service and networks, bio-inspired modelling of distributed network attacks was studied by \cite{Salgu}, and \cite{Chifu} optimized the semantic web service by bio-inspired methods.

Multiobjective optimization has also been addressed. For example,
\cite{Zeng2013} used non-dominated sorting genetic algorithm to solve constrained optimisation problems, while \cite{Yang2011BA} presents a multiobjective extension to the
bat algorithm.

Obviously, there are many other application topics, and interested readers can look at the table of contents of IJBIC for details.

\section{Recent Trends}

The rapid developments in bio-inspired computation mean that it is difficult to
predict its trends, though the current research activities may provide, to a certain extent, a good indication of what may happen next.

There is no single paper that can represent the current trends, though review papers tend to summarize more comprehensively what topics and applications have been addressed in recent studies. In this sense, review papers may provide a more comprehensive starting point.
For example, \cite{Akerhar} provided a relatively comprehensive review on bio-inspired computing, while \cite{Yang2011rev} provided an insightful analysis of randomization techniques in nature-inspired algorithms. In addition,
\cite{Parpi} surveyed new inspirations in swarm intelligence,
and \cite{Yang2013BAT} reviewed the bat algorithm comprehensively with detailed literature about various applications.

So what are the trends (you may wonder)? In fact, this is a very difficult question to answer, and we do not wish to mislead. Therefore, our comments below act as an optional,
personal notes, rather than an answer to this challenging questions. From our observations and the analysis of the table of contents of IJBIC, we may summarize the current trends as follows:
\begin{itemize}
\item Complex, real-world applications. More and more studies have focussed on real-world applications such as highly nonlinear design problems in engineering and industry. These applications tends to be complex with diverse and stringent constraints, and thus the problems can typically be multimodal.

\item Data intensive applications: the data volumes are increasing dramatically, driven by the information technology and social networks. Thus, data intensive data mining
    techniques tend to combine with bio-inspired algorithms such as PSO, cuckoo search and firefly algorithm to carry out fault detection and filtering as well as image processing.

\item Computationally expensive methods: Even with the best optimisation algorithms, the computational costs are usually caused by the high expense of evaluating objective functions, often in terms of external finite element or finite volume solvers,
    as those in aerospace and electromagnetic engineering. Therefore, approximate methods to save computational costs in function evaluations are needed.

\item Network and systems: Many research activities also focused on applications to complex networks and systems, including computer networks, electricity/energy networks, supply chain networks, and biological/ecological systems as well as social networks.

\item Novel applications: Sometimes, the existing methods and also new algorithms can be applied to study very new problems. Such novel applications can help to solve interesting and real-world problems. For example, bio-inspired algorithms can often be combined with traditional approaches to carry out classifications, feature selection
    and even combinatorial optimization such as the travelling salesman problem.

\end{itemize}
Obviously, as we mentioned earlier, there are other trends as well. Developments of new hybrid algorithms and continuing improvements of existing algorithms form a constant trend in the last twenty years. The above few points are just the parts that we hope researchers will continue to spend more time on, because these applications are important topics that will have a huge  impact in real-world applications.

\section{Challenges and Key Issues}

Despite the success in bio-inspired computation, many challenging problems remain unsolved.
Though different researchers may think differently in terms of the challenging problems,
however, we will highlight a few areas that further research and progress will be
more likely to have a huge impact.

\begin{itemize}
\item Theoretical analysis: Despite the rapid developments in applications, theoretical analysis still lacks behind. Any theoretical analysis will help to gain insight into the working mechanisms of bio-inspired algorithms, and certainly results concerning stability and convergence are useful. However, many algorithms still remain to be analyzed by mathematical tools. Future research efforts should focus more on the theoretical analysis.

\item Large scale problems: The applications in IJBIC are diverse, however, the vast majority of the application papers have dealt with small or moderate scale problems. The number of
    design variables are typically less than a hundred. As real-world applications are large scale with thousands or even millions of design variables, researchers should address complex large-scale  problems to test the scalability of optimisation algorithms and to produce truly useful results.

\item Combinatorial problems: Among complex optimisation problems, combinatorial optimization is typically hard to solve. In fact, most combinatorial problems do not usually have good methods to tackle with. In many cases, nature-inspired metaheurstic algorithms such as ant colony optimisation (ACO) and cuckoo search can be a good alternative. In fact, ACO, firefly algorithm and cuckoo search have all been applied to study the travelling salesman problem with good results. In reality, combinatorial problems can be from very different areas and applications, more research is highly needed.

\item True intelligence in algorithms: As the significant developments expand in
bio-inspired computation, some researchers may refer to some algorithms as `intelligent algorithm'. However, care should be taken when interpreting the meaning. Despite the names, algorithms are not truly intelligent at the moment. Some algorithms with fine-tuned performance and integration with expert systems may seem to show some low-level, basic intelligence, but they are far from truly intelligent. In fact, truly intelligent algorithms are yet to be developed.

\end{itemize}

Obviously, there are other important open problems, and some start to emerge with
interesting ideas. For example, all algorithms have some parameters, and the tuning of these parameters is important to the performance of an algorithm. There is no easy way to
tune an algorithm properly. A novel framework for self-tuning algorithms has been proposed by \cite{YangSTA}, which can tune an algorithm automatically.

\section{Concluding Remarks}

The rapid developments of bio-inspired computation can be reflected in the progress of the
Inderscience journal: IJBIC in the last five year. Despite the short history, IJBIC has
achieved successfully its main aim as a platform for exchanging research and ideas in bio-inspired computation. The diverse applications indicated the activeness and hotness of the topics, while the rarity in some topics such as theoretical results also poses further challenges for the journal.

From our observations and experience, we suggest the follow topics for further research and
welcome more papers in IJBIC:
\begin{itemize}

\item Theoretical analysis of bio-inspired algorithms.

\item Large-scale, real-world applications.

\item Combinatorial optimisation.

\item Data mining and image processing.

\item Novel applications.

\end{itemize}

This review paper celebrates the five year achievements of this journal and also highlights the future challenges and key issues. There is no doubt that IJBIC will continue its success with greater achievements in the coming years. It is hoped that this paper can inspire more research in both theory and applications in the future.


\begin{thebibliography}{100.}

\def\IJBIC{{\it Int. J. Bio-Inspired Computation}}
\def\vv#1#2{vol. #1, no. #2}

\bibitem[Adamatzky(2012)]{Adamatzky}
Adamatzky, A., (2012). `Slime mould computes planar shapes',
\IJBIC, \vv{4}{3}, pp. 149--154.

\bibitem[Akerkar and Sajja(2009)]{Akerhar}
Akerkar, R. and Sajja, P. S., (2009). `Bio-inspired computing: consituents and challenges',
\IJBIC, \vv{1}{3}, pp. 135--150.

\bibitem[Arumugam et al.(2009)]{Arumugam}
Arumugam, M. S., Murthy, G. R., C. K. Loo, (2009). `On the optimal control of the steel annealing processes as a two-stage hybrid systems via PSO algorithms',
\IJBIC, \vv{1}{3}, pp. 198--209.

\bibitem[Bessedik et al.(2011)]{Bessedik}
Bessedik, M., Toufik, B., Drias, H., (2011). `How can bees colour graphs',
\IJBIC, \vv{3}{1}, pp. 67--76.


\bibitem[Bonyadi and Shah-Hosseini(2010)]{Bonyadi}
Bonyadi, M. R, and Shah-Hosseini, H., (2010). `A dynamic max-min ant system for solving
the travelling salesman problem', \IJBIC, \vv{2}{6}, pp. 422-433.


\bibitem[Cai et al.(2012)]{Cai2012}
Cai, X.J., Fan, S., and Tan, Y., (2012). `Light responsive curve selection for photosynthesis
operator of APOA', \IJBIC, \vv{4}{6}, pp. 373--379.


\bibitem[Chifu et al.(2013)]{Chifu}
Chifu, V. R., Pop, C. B., Solomie, I., Suia, D. S., Niculici, A., Negrean, A., Jeflea, H., (2013).
`Optimising the semantic web service composition process using bio-inspired methods',
\IJBIC, \vv{5}{4}, pp. 226--238.

\bibitem[Cui and Cai(2009)]{Cui}
Cui Z. H., and Cai X. J. (2009)  `Integral particle swarm optimisation with dispersed
accelerator information', {\it Fundam. Inform.}, \vv{95}{2}, 427--447.

\bibitem[Cui et al.(2010)]{Cui2010}
Cui, Z. H., Cai, X. J., Zeng, J. C., Yin, Y., (2010). `PID-controlled particle swarm optimization',
{\it Multiple-Valued Logic and Soft Computing}, \vv{16}{6}, pp. 585--609.


\bibitem[Cui et al.(2013)]{Cui2013}
Cui, Z. H., Fan, S. J., Zeng, J. C., Shi, Z.,Z., (2013). `Artificial plan optimisation
algorithm with three-period photosynthesis', \IJBIC, \vv{5}{2}, pp. 133--139.

\bibitem[Deb(2009)]{Deb2009}
Deb, S., `A novel alternative to conventional machince vision features',
\IJBIC, \vv{1}{1/2}, pp.89--92.

\bibitem[Dey et al.(2013)]{Dey}
Dey, N., Samanta, S., Yang, X. S., Das, A., Chaudhuri, S. S., (2013). `Optimisation of scaling factors in
electrocardiogram signal watermarking using cuckoo search', \IJBIC, \vv{5}{5}, pp. 315--326.

\bibitem[Fister et al.(2013a)]{Fister}
Fister, I., Fister Jr., I., Yang, X. S., Brest, J., (2013).
`A comprehensive review of firefly algorithms', {\it Swarm and Evolutionary Computation},
\vv{13}{1}, pp. 34--46.


\bibitem[Fister et al.(2013b)]{Fister2013b}
Fister, I., Yang, X. S., Brest, J., Fister Jr., I., (2013b). `Modified firefly algorithm using
quaternion representation', {\it Expert Syst. Appl.}, \vv{40}{18}, pp. 7220-7230.

\bibitem[Gandomi et al.(2013a)]{Gandomi2013CPSO}
Gandomi, A. H., Yun, G. J., Yang, X. S., Talatahari, S. (2013a).
`Chaos-enhanced accelerated particle swarm optimization', {\it
Communications in Nonlinear Science and Numerical Simulation}, \vv{18}{2},
pp. 327--340 (2013).

\bibitem[Gandomi et al.(2013b)]{Gandomi2013FA}
Gandomi, A. H., Yang, X. S., Talatahari, S., Alavi, A. H., (2013b).
`Firefly algorithm with chaos', {\it Communications in Nonlinear Science and Numerical Simulation},
\vv{18}{1}, pp. 89--98.

\bibitem[Green(2011)]{Green}
Green, D. G., `Elements of a network theory of complex adaptive systems',
\IJBIC, \vv{3}{3}, pp. 159--167.

\bibitem[Gross and Dorigo(2009)]{GrossDorigo}
Gross, R. and Dorigo, M., (2009). `Towards group transport by swarms of robots',
\IJBIC, \vv{1}{1/2}, pp. 1--13.

\bibitem[Handa(2012)]{Handa}
Handa, H., (2012). `Neuroevolution with manifuold learning for playing Mario',
\IJBIC, \vv{4}{1}, pp. 14--26.


\bibitem[Jamil and Zepernick(2013)]{Jamil}
Jamil, M., Zepernick, H. J., (2013). `Multimodal function optimisation
with cuckoo search algorithm', \IJBIC, \vv{5}{2}, pp. 73--78.

\bibitem[Kennedy and Eberhart(1995)]{Ken}
Kennedy J. and Eberhart R. C., (1995). `Particle swarm optimization'. {\it
Proc. of IEEE International Conference on Neural Networks},
Piscataway, NJ. pp. 1942-1948.

\bibitem[Khabbazi et al.(2009)]{Khabb}
Khabbazi, A., Gargari, E. A., Lucas, C., (2009). `Imperialist competitive algorithm for minimum bit
error rate beamforming', \IJBIC, \vv{1}{1/2}, pp. 125--133.

\bibitem[Komatsu and Namatame(2011)]{Komatsu}
Komatsu, T., Namatame, A. (2011). `Dynamic diffusion in evolutionary optimised networks',
\IJBIC, \vv{3}{6}, pp. 384--392.

\bibitem[Laalaoui and Drias(2010)]{Laalaoui}
Laalaoui, Y. and Drias, H., (2010). `ACO approach with learning for preemptive scheduling of
real-time tasks', \IJBIC, \vv{2}{6}, pp. 383--394.

\bibitem[Layeb(2011)]{Layeb}
Layeb, A., (2011). `A novel quantum inspired cuckoo search for knapsack problems',
\IJBIC, \vv{3}{5}, pp. 297--305.

\bibitem[Liu and Zhang(2012)]{Liu2012}
Liu, J. G., Fang, N., (2012). `Computational intelligence and its applications in engineering: editorial', \IJBIC, \vv{4}{2}, pp. 61--62.

\bibitem[Marichelvam(2012)]{Marichel}
Marichelvam, M. K., (2012). `An improved hybrid cuckoo search metaheuristics algorithm for permutation
flow shop scheduling problems', \IJBIC, \vv{4}{4}, pp. 200--205.

\bibitem[Natarajan et al.(2012)]{Natar}
Natarajan, A., Subramanian, S., Prematatha, K., (2012). `A comparative study of cuckoo search
and bat algorithm for bloom filter optimisation in spam filtering', \IJBIC, \vv{4}{2}, pp. 89--99.

\bibitem[Neri and Caponio(2010)]{Neri}
Neri, F., Caponio, A., (2010). `A differential evolution for optimisation in noisy environment',
\IJBIC, \vv{2}{3/4}, pp. 152--168.

\bibitem[Pal et al.(2011)]{Pal2011}
Pal, S., Basak, A., Das, S., (2011). `Linear antenna array synthesis with modified invasive
weed optimisation algorithm', \IJBIC, \vv{3}{4}, pp. 238--251.

\bibitem[Pandi et al.(2010)]{Pandi}
Pandi, V. R., Panigrahi, B. K., Das, S., Cui, Z. H., (2010).
`Dynamic economic load dispatch with wind energy using modified harmony search',
\IJBIC, \vv{2}{3/4}, pp.282--289.


\bibitem[Parpinelli and Lopes(2011)]{Parpi}
Parpinelli, R. S. and Lopes, H. S., (2011).
`New inspirations in swarm intelligence: a survey', \IJBIC, \vv{3}{1}, pp. 1--16.

\bibitem[Ray and Mondal(2011)]{Ray2011}
Ray, K. S, and Mondal, M., (2011). `Similarity-based fuzzy reasoning by DNA computing',
\IJBIC, \vv{3}{2}, pp. 112--122.

\bibitem[Roy et al.(2010)]{Roy2010}
Roy, G. G., Chakraborty, P., Das, S., `Designing fraction-order PI$^{\lambda}$D$^{\mu}$ controller
using differential harmony search algorithm', \IJBIC, \vv{2}{5}, pp. 303--309.

\bibitem[Saha et al.(2013)]{Saha}
Saha, S. K., Mandal, D., Ghoshal, S. P., (2013). `Bacteria foraging optimisation algorithm
for optimal FIR filter design', \IJBIC, \vv{5}{1}, pp. 52--66.

\bibitem[Solgueiro and Abreu(2013)]{Salgu}
Salgueiro, P., Abreu, S., (2013). `Modelling distributed network attacks with constraints',
\IJBIC, \vv{5}{4}, pp. 210--225.

\bibitem[Shah-Hosseini(2009)]{Shah}
Shah-Hosseini, H., (2009). `The intelligence water drops algorithm:
a nature-inspired swarm-based optimization algorithm, \IJBIC, \vv{1}{1/2}, pp. 71-79.

\bibitem[Sheta(2010)]{Sheta}
Sheta, A. F., (2010). `Analogue filter design using differential evolution',
\IJBIC, \vv{2}{3/4}, pp. 233-241.

\bibitem[Sivanandam and Visalakshi(2009)]{Sivan}
Sivanandam, S. N., Visalakshi, P., (2009).
`Dynamic task scheduling with load balancing using parallel
orthogonal particle swarm optimisation', \IJBIC, \vv{1}{4}, pp. 276--286.

\bibitem[Singh et al.(2010)]{Singh}
Singh, N. A., Muraleedharan, K. A., Gomathy, K., (2010).
`Damping of low frequency oscillations in power system network using swarm intelligence
tuned fuzzy controller',\IJBIC, \vv{2}{1}, pp. 1--8.

\bibitem[Sing et al.(2013)]{Sing2013}
Sing, M., Ketan, B., Abhyankar, A. R., (2013). `Optimal coordination of electro-mechanical-based
overcurrent relays using artificial bee colony algorithm', \IJBIC, \vv{5}{5}, pp. 267--280.

\bibitem[Srivastava et al.(2012a)]{Srivastava}
Srivastava, P. R., Sravya, C., Kamisett, S., Lakshmi, M., (2012a).
`Test sequence optimisation: an intelligent approach via cuckoos search',
\IJBIC, \vv{4}{3}, pp. 139--148.

\bibitem[Srivastava et al.(2012b)]{Srivastava2012b}
Srivastava, P. R., Varashney, A., Nama, P, Yang, X. S., (2012b).
`Software test effort estimation: a model based on cuckoo search', \IJBIC, \vv{4}{5}, pp. 278--285.

\bibitem[Srivastava et al.(2013)]{SrivSoft}
Srivastava, P. R., Mallikarjun, B., Yang, X. S., (2013). `Optimal test sequence generation using firefly
algorithm', {\it Swarm and Evolutionary Computation}, \vv{8}{1}, pp. 44--53.


\bibitem[Tanimoto(2011)]{Tanimoto}
Tanimoto, J., (2011). 'Evolving world: editorial', \IJBIC, \vv{3}{3}, pp. 141-141.


\bibitem[Xiao and Chen(2013)]{Xiao}
Xiao, R. B. and Chen, T. G., (2013). `Relatioships of swarm intelligence and artificial
immune system', \IJBIC, \vv{5}{1}, pp. 35--51.

\bibitem[Xie et al.(2010)]{Xie2010}
Xie, L. P., Tan, Y., Zeng, J. C., Cui, Z. H., (2010).
`Artificial physics optimisation: a brief survey', \IJBIC, \vv{2}{5}, pp. 291--302.


\bibitem[Yampolskiy(2010)]{Yampol}
Yampolskiy, R. V., (2010). `Application of bio-inspired algorithm to the problme of integer
factorisation', \IJBIC, \vv{2}{2}, pp. 115--123.

\bibitem[Yang(2008)]{Yang}
Yang, X. S., (2008). {\it Nature-Inspired Metaheuristic Algorithms}, 1st Edition,
Luniver Press, Frome, UK.


\bibitem[Yang(2010)]{Yang2010FA}
Yang, X. S., (2010). Firefly algorithm, stochastic test functions and design optimisation,
\IJBIC, \vv{2}{2}, pp. 78--84.

\bibitem[Yang et al.(2013)]{YangSTA}
Yang, X. S., Deb, S., Loomes, M., Karamanoglu, M., (2013).
`A framework for self-tuning optimization algorithm', {\it Neural Computing and Applications},
\vv{23}{7/8}, pp. 2051--2057.

\bibitem[Yang(2013)]{Yang2013FA}
Yang, X. S., `Multiobjective firefly algorithm for continuous optimization',
{\it Engineering with Computers}, \vv{29}{2}, pp. 175--184.

\bibitem[Yang(2011a)]{Yang2011rev}
Yang, X. S., (2011a). `Review of metaheuristics and generalized evolutionary walk algorithm',
{\it Int. J. Bio-Inspired Compuation}, \vv{3}{2}, pp. 77--84.

\bibitem[Yang(2011b)]{Yang2011BA}
Yang, X. S., (2011b). `Bat algorithm for multi-objective optimisation', \IJBIC, \vv{3}{5},
pp. 267--274.

\bibitem[Yang(2012)]{Yang2012SI}
Yang, X. S., (2012). `Metaheuristics and swarm intelligence in engineering
and industry: editorial', \IJBIC, \vv{4}{4}, pp. 197--199.

\bibitem[Yang and Deb(2012)]{YangES}
Yang, X. S. and Deb, S., (2012). `Two-stage eagle strategy with differential evolution',
\IJBIC, \vv{4}{1}, pp. 1--5.

\bibitem[Yang and He(2013)]{Yang2013BAT}
Yang, X. S., He, S., (2013). `Bat algorithm: literature review and applications',
\IJBIC, \vv{5}{3}, pp. 141--149.

\bibitem[Yang and Koziel(2010)]{YangKoziel}
Yang, X. S. and Koziel, S., (2010).
`Computational optimization, modelling and simulation -- a paradigm
shift', {\it Procedia in Computer Science}, \vv{1}{1}, pp. 1291--1294.

\bibitem[Yang et al.(2013)]{Yangbook2013}
Yang, X. S., Cui, Z. H., Xiao, R. B., Gandomi, A. H., Karamanoglu, M., (2013).
{\it Swarm Intelligence and Bio-Inspired Computation: Theory and Applications},
Elsevier, London.

\bibitem[Zeng et al.(2013)]{Zeng2013}
Zeng, S. Y., Zhou, D., Li, H., (2013). `Non-dominated sorting genetic algorithm with
decomposition to solve constrained optimisation problems', \IJBIC, \vv{5}{3}, pp. 150--163.

\bibitem[Zhang and Cai(2010)]{ZhangCai}
Zhang, T. and Cai, J. D., (2010). A novel hybird particle swarm optimisation method applied to
economic dispatch, \IJBIC, \vv{2}{1}, pp. 9--17.

\end{thebibliography}
\end{document}